\documentclass[12pt,leqno]{article}

\newcounter{conjecture}
\newcounter{remark}
\newcounter{corollary}

\newcommand{\eqnsection}{
    \renewcommand{\theequation}{\thesection.\arabic{equation}}
    \makeatletter
    \csname @addtoreset\endcsname{equation}{section}
    \makeatother}
\newtheorem{theorem}{Theorem}
\newtheorem{lemma}{Lemma}

\newcommand{\dd}{\delta}

\newcommand{\lar}{\longrightarrow}
\newcommand{\eps}{\varepsilon}
\newcommand{\om}{\Omega}
\newcommand{\aaa}{\alpha}
\newcommand{\reals}{R}
\newcommand{\lll}{\label}

\def \be{\begin{equation}}
\def \ee{\end{equation}}
\def \bt{\begin{theorem}}
\def \et{\end{theorem}}
\def \bea{\begin{eqnarray}}
\def \eea{\end{eqnarray}}
\def \bas{\begin{eqnarray*}}
\def \eas{\end{eqnarray*}}

\def \al{\alpha}

\def \om{\omega}

\def \({\left(}
\def \){\right)}

\def \nn{\nonumber}

\def \bc{\begin{center} }
\def \ec{\end{center} }
\def \bs{\begin{slide} }
\def \es{\end{slide} }

\def\square{{\vcenter{\vbox{\hrule height.3pt
         \hbox{\vrule width.3pt height5pt \kern5pt
            \vrule width.3pt}
         \hrule height.3pt}}}}
\def\qed{{\hfill $\square$ \bigskip}}

\eqnsection
\begin{document}

\title{A Tanaka formula for the derivative of intersection local time in $\reals^1$}

\author{Greg Markowsky}

\maketitle

\begin{abstract}
Let $B_t$ be a one dimensional Brownian motion, and let $\aaa'$
denote the derivative of the intersection local time of $B_t$ as
defined in \cite{rosen2}. The object of this paper is to prove the
following formula
\be \nn \frac{1}{2}\aaa_{t}'(x) + \frac{1}{2}sgn(x)t =
 \int_0^t L_s^{B_s - x}dB_s - \int_{0}^{t} sgn
(B_{t}-B_{u} - x) du \ee
which was given as a formal identity in \cite{rosen2} without proof.

\end{abstract}

Let $B$ denote Brownian motion in $R^1$. In \cite{rosen2}, Rosen
demonstrated the existence of of a process which he termed the
derivative of self intersection local time for $B$. That is, he
showed that there is a process $\alpha_t(y)$, formally defined as

\be \alpha_t(y) = -\int_0^t \int_0^s \dd '(B_s -B_r -y)drds \ee such
that, for any $C^1$ function $g$, we have

\be \int_0^t \int_0^s g'(B_s -B_r -y)drds = - \int_R g(y) \al '_t(y)
dy \ee In this paper we'll prove a Tanaka-style formula for
$\alpha'$ which was given without proof by Rosen in \cite{rosen2}.
We define

\be sgn(x) = \left \{ \begin{array}{ll}
-1 & \qquad  \mbox{if } x < 0  \\
0 & \qquad \mbox{if } x = 0 \\ 1 & \qquad \mbox{if } x > 0
\end{array} \right. \ee Our result is

\bt There is a set of measure one upon which the following holds for
all x and t:
\be \frac{1}{2}\aaa_{t}'(x) + \frac{1}{2}sgn(x)t =
 \int_0^t L_s^{B_s - x}dB_s - \int_{0}^{t} sgn
(B_{t}-B_{u} - x) du \ee
 \et
\medskip
{\bf Proof:} Fix $t$ and $x$ for the time being. In what follows,
the constant $c$ may change from line to line. Let $f(x) = \pi
^{-1/2} e^{-x^2}$. Let $f_{\eps}(x) = \frac{1}{\eps}
f(\frac{x}{\eps})$, so that $f_{\eps} \lar \dd$ weakly as $\eps \lar
0$. We assume in all calculations below that $\eps < 1$. Let \bea
F_{\eps}(x) = \int_{0}^{x} f_{\eps}(t) dt &=&
\int_{0}^{\frac{x}{\eps}}f(t) dt \eea We apply Ito's formula to
$F_{\eps}$ to get \bea && F_{\eps}(B_{t} - B_{u} - x) - F_{\eps}(-x)
= \int_{u}^{t}f_{\eps}(B_{s} - B_{u} - x) dB_{s} + \nn\\&&
\frac{1}{2} \int_{u}^{t} f_{\eps}'(B_{s} - B_{u} - x) ds \eea which
gives \bea \lll{emo} && \int_{0}^{t}F_{\eps}(B_{t} - B_{u} - x)du -
t F_{\eps}(-x) \nn\\&& = \int_{0}^{t}\int_{0}^{s}f_{\eps}(B_{s} -
B_{u} - x) du dB_{s} + \frac{1}{2} \int_{0}^{t} \int_{u}^{t}
f_{\eps}'(B_{s} - B_{u} - x) ds du \eea Note that $F_{\eps}(x) \lar
\frac{1}{2}sgn(x)$ as $\eps \lar 0$. Furthermore, $|F_{\eps}(x)|
\leq \frac{1}{2}$ for all $x, \eps$, so by the dominated convergence
theorem, the first integral on the left approaches $\int_{0}^{t}
sgn(B_{t} - B_{u} - x) du $ as $\eps \lar 0$. By Theorem 1 in
\cite{rosen2}, the rightmost integral on the right side is equal to

\be -\frac{1}{2} \int_{R} f_{\eps}(y-x) \aaa_{t}'(y)dy \ee This term
approaches $-\frac{1}{2}\aaa_{t}'(x)$ as $\eps \lar 0$ for all $x$
at which $\aaa_{t}'(x)$ is continuous. In \cite{rosen2} it was shown
that $\aaa_{t}'(x)$ is continuous for all $x \neq 0$. To deal with
the case $x=0$, we need another fact proved in \cite{rosen2}, namely
that $\aaa_{t}'(x) + sgn(x)$ is continuous in $x$. Using this,
together with the fact that $f_{\eps}(x) sgn(x)$ is an odd function,
we have the following string of equalities: \bea \lim_{\eps \lar 0}
\int_{R} f_{\eps}(y) \aaa_{t}'(y)dy \\ \nn = \lim_{\eps \lar 0}
\int_{R} f_{\eps}(y) (\aaa_{t}'(y)+sgn(y))dy
\\ \nn = \aaa_{t}'(0)+sgn(0) = \aaa_{t}'(0) \eea The only term which remains is the
leftmost term on the right side of (\ref{emo}):

\be V(x,\eps) := \int_{0}^{t}\int_{0}^{s}f_{\eps}(B_{s} - B_{u} - x)
du dB_{s} \ee We will show that

\be \lll{yar} \int_{0}^{s}f_{\eps}(B_{s} - B_{u} - x) du \lar
L_s^{B_s - x} \ee in $L^2$, and this is enough to complete the proof
for fixed $x$ and $t$. Using the standard occupation times formula,
we have a.s.

\be \int_{0}^{s}f_{\eps}(B_{s} - B_{u} - x) du = \int f_{\eps}(B_s -
y - x) L^y_s dy \ee Since $\int f_{\eps} = 1$, we have
\bea \lll{davai} E[\int f_{\eps}(B_s - y - x)L^y_s dy - L_s^{B_s - x}]^2 \\
\nn \leq E[\int f_{\eps}(B_s - y - x)|L^y_s - L_s^{B_s - x}| dy]^2
\\ \nn \leq E[\int f_{\eps}(B_s - y - x)|L^y_s - L_s^{B_s - x}|^2dy]
\eea The last inequality is Jensen's inequality, as $f_{\eps}(B_s -
y - x)dy$ is a probability measure on $\reals$. We integrate
separately over the two regions $\{ |y - (B_s - x)| < \sqrt{\eps}
\}$ and $\{ |y - (B_s - x)| \geq \sqrt{\eps} \}$. We can bound the
contribution to (\ref{davai}) from the second region by

\be 2 E[\int_{\{ |y - (B_s - x)| \geq \sqrt{\eps \}}} f_{\eps}(B_s -
y - x)(|L^y_s|^2 + |L_s^{B_s - x}|^2)dy] \ee Expand this into the
expectation of two integrals. The first is

\be \lll{sw} E[\int_{\{ |y - (B_s - x)| \geq \sqrt{\eps \}}}
f_{\eps}(B_s - y - x)(L^y_s)^2 dy] \ee Since $ |y - (B_s - x)| \geq
\sqrt{\eps}$, we see that $f_{\eps}(B_s - y - x) \leq c (1/\eps)
e^{-1/\eps}$. Thus, (\ref{sw}) is bounded by

\be \lll{hil} c (1/\eps) e^{-1/\eps} \int E[(L^y_s)^2] dy \leq  c
(1/\eps) e^{-1/\eps} \ee We have used here the fact that $\int
E[(L^y_s)^2] dy < \infty$. One way of proving this is to note that
$E[(L^y_s)^2] \leq P(T_y < s) E[(L^0_s)^2]$ by the strong Markov
property, where $T_y$ is the first hitting time of $y$. $P(T_y < s)
= P[|B_s| > |y|]$ by the reflection principle, and it is
straightforward to check that
\be \int P[|B_s| > |y|] dy < \infty \ee
Thus, (\ref{hil}) converges to $0$ as $\eps \lar 0$. The second
integral is \bea \lll{nir} E[ |L_s^{B_s - x}|^2 \int_{\{ |y - (B_s -
x)| \geq \sqrt{\eps \}}} f_{\eps}(B_s - y - x)dy] \\ \nn \leq
\int_{|y|>\sqrt{\eps}} f_{\eps}(y)dy E[|L_s^{B_s - x}|^2 ] \eea We
require the fact that $E[|L_s^{B_s - x}|^2 ]$ is finite, and this
may be proved as follows:
\bea \lll{filt} L_s^{B_s - x} = \lim_{\eps \lar 0} \int_0^s f_{\eps}((B_s-B_u) - x) du = \\
\nn \lim_{\eps \lar 0}\int_0^s f_{\eps}((B_s-B_{s-u}) - x) du =
\tilde{L}_s^x \eea where $\tilde{L}_s^x$ is the local time of the
Brownian motion $\tilde{B}_u = (B_s-B_{s-u})$. Then
$E[(\tilde{L}_s^x)^2]$ is bounded by $E[(\tilde{L}_s^0)^2]< \infty$,
for \be \lll{offs} \tilde{L}_s^x =_{law} 1_{\{\tilde{T}_x < s\} }
\tilde{L'}_{s-\tilde{T}_x}^0\ee where $\tilde{T}_x$ is the first
time $\tilde{B_u}$ hits $x$, and $\tilde{L'}$ is the local time of
the Brownian motion $\tilde{B'}_u = \tilde{B}_{T_x +u} -
\tilde{B}_{T_x}$. (\ref{offs}) is a.s. smaller than $\tilde{L}_s^0$,
as local time is increasing in $s$. Thus, (\ref{nir}) is bounded by

\be c \int_{|x|>\sqrt{\eps}} f_{\eps}(x)dx = c
\int_{|x|>\eps^{-1/2}} f(x)dx \ee and this approaches $0$ as $\eps
\lar 0$. We must now show that

\be \lll{vel} E[\int_{|B_s - y - x| < \sqrt{\eps}} f_{\eps}(B_s - y
- x)|L^y_s - L_s^{B_s - x}|^2dy] \ee approaches $0$ as $\eps$ does.
We will need the following lemma.

\begin{lemma} Given $\dd >0$, there is an $M>0$ such that

\be E[(L_s^{B_s-x})^2 1_{\{|B_s - x| > M \} }] < \dd \ee

\end{lemma}

{\bf Proof:} By the Cauchy-Schwarz inequality, we have \be
E[(L_s^{B_s-x})^2 1_{\{|B_s - x| > M \} }] \leq
E[(L_s^{B_s-x})^4]^{1/2}  P(|B_s - x| > M)^{1/2} \ee Writing
$\tilde{L}_s^x$ for $L_s^{B_s-x}$ as we have done before, we see \be
E[(L_s^{B_s-x})^4] = E[(\tilde{L}_s^{x})^4] \leq
E[(\tilde{L}_s^{0})^4]\ee with the last inequality being due to the
same argument as in steps (\ref{filt}) and (\ref{offs}). Local time
at $0$ has all moments, so $E[(L_s^{B_s-x})^4]$ is uniformly
bounded. It is evident that $P(|B_s - x| > M) \lar 0$ as $M \lar
\infty$. This proves the lemma. \qed

Fix $\dd > 0$. The lemma, together with the fact that \be
E[(L_s^y)^2] \leq P(T_M < s) E[(L_s^0)^2] \ee when $y > M$, allows
us to pick $M$ sufficiently large so that

\be E[(L_s^y)^2], E[(L_s^{B_s-x})^2 1_{\{|B_s - x| > M \} }] < \dd
\ee when $y>M$. Then, substituting $y' = y - (B_s - x)$ \bea
\lll{vel2} E[\int_{\{|B_s - y - x| < \sqrt{\eps}\} \bigcap \{
|y|>M+1\} } f_{\eps}(B_s - y - x)|L^y_s - L_s^{B_s - x}|^2dy] \\
\nn = E[\int_{\{|y'| < \sqrt{\eps}\} \bigcap \{ |y' + (B_s
-x)|>M+1\} } f_{\eps}(y')|L^{y' +(B_s - x)}_s - L_s^{B_s -
x}|^2dy']\\ \nn \leq c\int_R f_{\eps}(y')E[|L^{y' +(B_s -
x)}_s|^21_{\{ |y' +
(B_s -x)|>M+1\}} + |L_s^{B_s - x}|^21_{\{|B_s - x| > M \} }]dy' \\
\nn \leq \dd c \int_{\reals} f_{\eps}(y') dy' = \dd c \eea Therefore
we can restrict the integral to the region $|y|<M+1$, which means
$|B_s - x| < M + 2$.  Now, by \cite{rosmar} there is an $L^2$ random
variable $X(\om)$ such that $|L_s^y - L_s^z| \leq X(\om)|y - z|^k$,
where $k>0$ is any number less than $1/2$, whenever $|y|, |z| <
M+2$. Using this we have \bea \lll{vel3} E[\int_{\{|B_s - y - x| <
\sqrt{\eps}\} \bigcap \{
|y|<M+1\} } f_{\eps}(B_s - y - x)|L^y_s - L_s^{B_s - x}|^2dy] \\
\nn \leq \eps^{k/2} E[X(\om)^2 \int_{\reals} f_{\eps}(B_s - y -
x)dy]\eea The $dy$ integral is bounded by $1$, so (\ref{vel3}) is
bounded by $\eps^{2k} E[X(\om)^2]$, and this converges to $0$ as
$\eps$ goes to $0$. This proves that

\be \lll{vel2} E[\int_{|B_s - y - x| < \sqrt{\eps}} f_{\eps}(B_s - y
- x)|L^y_s - L_s^{B_s - x}|^2dy] \lar 0\ee as $\eps \lar 0$, and
proves the result for fixed $x,t$. We would like to prove it to be
true for all $x,t$ on a set of full measure, however. We will do so
by proving that, for $t,t' < M$ we have

\be \lll{req} E[V(x,\eps,t)-V(x',\eps ',t')]^{2n} \leq C_M
|(x,\eps,t)-(x',\eps',t')|^{n/20}\ee for any positive integer $n
\geq 3$. This will allow us to apply Kolmogorov's criteria(see
\cite{revyor}, Theorem I.2.1) for uniform continuity, which will
complete the proof. We will in fact show separately that

\be \lll{req2} E[V(x,\eps,t)-V(x',\eps ,t)]^{2n} \leq C_M
|x-x'|^{2n/3}\ee

\be \lll{req3} E[V(x,\eps,t)-V(x,\eps' ,t)]^{2n} \leq C_M
|\eps-\eps'|^{2n/3}\ee

\be \lll{req4} E[V(x,\eps,t)-V(x,\eps ,t')]^{2n} \leq C_M
|t-t'|^{(n-1)/10}\ee
and these clearly imply (\ref{req}). In order
to prove this, we'll need a convenient expression bounding
$E(V(x,\eps,t))^{2n}$. We'll use the identity

\be f_{\eps}(x) = \frac{i}{2\pi} \int_{\reals} e^{ixp} \hat{f}(\eps
p)dp \ee
By the Burkholder-Davis-Gundy inequality(again see
\cite{revyor}, Corollary IV.4.2) we have
\bea && E(V(x,\eps,t))^{2n} \leq cE(\int_{0}^{t}(\int_{0}^{s}
\int_{\reals}e^{i(B_{s}-B_{u}-x)p} \hat{f}(\eps p)\, dp \, du)^{2}\,
ds)^{n} \nn\\ && = c\int_{\reals^{2n}}\int_{[0,t]^{n}}
\int_{[o,s]^{2n}}(\prod_{i=1}^n \hat{f}(\eps p_{i}) \hat{f}(\eps
p'_{i})) E[exp(i\sum_{i=1}^n
[p_{i}(B_{s_{i}}-B_{u_{i}}) + p'_{i}(B_{s_{i}}-B_{u'_{i}})] \nn\\
\nn && exp(ix\sum_{i=1}^n p_{i})exp(ix\sum_{i=1}^n
p'_{i})(\prod_{i=1}^n du_{i} du'_{i} ds_{i} dp_{i} dp'_{i}) \eea
where $i$ ranges from 1 to $n$ in the products
and sum. We will deal first with the variance in $x$ and $\eps$. We
have the following bounds:

\be |e^{ipx}-e^{ipx'}| \leq c|p|^{1/3}|x-x'|^{1/3} \ee
\be |\hat{f}(\eps p)-\hat{f}(\eps' p)| \leq
c|p|^{1/3}|\eps-\eps'|^{1/3} \ee We will also use the trivial bounds
$|e^{ipx}|, |\hat{f}(\eps p)| \leq 1$. Thus,
\bea  && E(V(x,\eps,t)-V(x',\eps,t))^{2n} \nn\\&& \leq
c\int_{\reals^{2n}}\int_{[0,t]^{n}} \int_{[o,s]^{2n}}(\prod_{i=1}^n
\hat{f}(\eps p_{i})\hat{f}(\eps p'_{i}))E[exp(i\sum_{i=1}^n
[p_{i}(B_{s_{i}}-B_{u_{i}}) + p'_{i}(B_{s_{i}}-B_{u'_{i}})] \nn\\
\nn && (\prod_{i=1}^n
|e^{ip_{i}x}-e^{ip_{i}x'}||e^{ip'_{i}x}-e^{ip'_{i}x'}|)(\prod_{i=1}^n
du_{i} du'_{i} ds_{i} dp_{i} dp'_{i}) \nn\\&& \leq
c|x-x'|^{2n/3}\int_{\reals^{2n}}\int_{[0,t]^{n}} \int_{[o,s]^{2n}}
E[exp(i\sum_{i=1}^n
[p_{i}(B_{s_{i}}-B_{u_{i}}) + p'_{i}(B_{s_{i}}-B_{u'_{i}})] \nn\\
\nn && (\prod_{i=1}^n |p_{i}|^{1/3} |p'_{i}|^{1/3})(\prod_{i=1}^n
du_{i} du'_{i} ds_{i} dp_{i} dp'_{i}) \nn \eea
Likewise,
\bea \lll{hoo} E(V(x,\eps,t)-V(x,\eps',t))^{2n}  \leq
c\int_{\reals^{2n}}\int_{[0,t]^{n}} \int_{[o,s]^{2n}}(\prod_{i=1}^n
exp(ix\sum_{i=1}^n p_{i})exp(ix\sum_{i=1}^n p'_{i})) \nn \\ \nn
E[exp(i\sum_{i=1}^n [p_{i}(B_{s_{i}}-B_{u_{i}}) +
p'_{i}(B_{s_{i}}-B_{u'_{i}})] (\prod_{i=1}^n |\hat{f}(\eps p_{i}) -
\hat{f}(\eps' p_{i})|)(\prod_{i=1}^n du_{i} du'_{i} ds_{i} dp_{i}
dp'_{i})\nn\\ \leq
c|\eps-\eps'|^{2n/3}\int_{\reals^{2n}}\int_{[0,t]^{n}}
\int_{[o,s]^{2n}} E[exp(i\sum_{i=1}^n
[p_{i}(B_{s_{i}}-B_{u_{i}}) + p'_{i}(B_{s_{i}}-B_{u'_{i}})] \nn\\
\nn (\prod_{i=1}^n |p_{i}|^{1/3} |p'_{i}|^{1/3})(\prod_{i=1}^n
du_{i} du'_{i} ds_{i} dp_{i} dp'_{i})\nn \eea
In order to control
the variance in $\eps$ and $x$ in the required (\ref{req}) we need
only bound \bea && \lll{wri} c\int_{\reals^{2n}}\int_{[0,t]^{n}}
\int_{[o,s]^{2n}} E[exp(i\sum_{i=1}^n
[p_{i}(B_{s_{i}}-B_{u_{i}}) + p'_{i}(B_{s_{i}}-B_{u'_{i}})] \\
\nn && (\prod_{i=1}^n |p_{i}|^{1/3} |p'_{i}|^{1/3})(\prod_{i=1}^n
du_{i} du'_{i} ds_{i} dp_{i} dp'_{i})\nn \eea The value of the
expectation in the integrand will depend on the ordering of the
$s_i$'s, $u_i$'s, and $u'_i$'s. For example, if $n=2$, then while
considering the region $s_{1}>s_{2}>u_{1}>u'_{1}>u_{2}>u'_{2}$, we
rewrite the integrand as
\bea \lll{ywa}
E[exp(i[(p_{1}+p'_{1})(B_{s_{1}}-B_{s_{2}}) +
(p_{1}+p'_{1}+p_{2}+p'_{2})(B_{s_{2}}-B_{u_{1}})+ \\
(p'_{1}+p_{2}+p'_{2})(B_{u_{1}}-B_{u'_{1}})+(p_{2}+p'_{2})(B_{u'_{1}}-B_{u_{2}})+
(p'_{2})(B_{u_{2}}-B_{u'_{2}})])] \nn \eea
By the independence of
increments of Brownian motion, this expectation splits, and is equal
to
\bea \lll{off} && exp(-[(p_{1}+p'_{1})^{2}(s_{1}-s_{2}) +
(p_{1}+p'_{1}+p_{2}+p'_{2})^{2}(s_{2}-u_{1})+ \nn \\&&
(p'_{1}+p_{2}+p'_{2})^{2}(u_{1}-u'_{1})+(p_{2}+p'_{2})^{2}(u'_{1}-u_{1})+
(p'_{2})^{2}(u_{1}-u'_{2})]) \eea We now substitute \bea &&
\lll{jpb} (z_1,z_2,z_3,z_4,z_5,z_6) =
(s_{1}-s_{2},s_{2}-u_{1},u_{1}-u'_{1},u'_{1}-u_{2},u_{2}-u'_{2},u'_2)\eea
and integrate with respect to the $z_i$'s using the simple bound
\be
\lll{bund} \int_{0}^{t}e^{-rb^{2}}dr \leq \frac{c}{1+b^{2}} \ee
We see that in order to show (\ref{wri}) is bounded in this case we
must show
\bea && \int_{\reals^{4}} \frac{(\prod_{i=1}^2 |p_{i}|^{1/3}
|p'_{i}|^{1/3})}{(1+(p_{1}+p'_{1})^{2})(1+(p_{1}+p'_{1}+p_{2}+p'_{2})^{2})(1+(p'_{1}+p_{2}+p'_{2})^{2})}
\nn\\&& \frac{(\prod_{i=1}^2
dp_{i}dp'_{i})}{(1+(p_{2}+p'_{2})^{2})(1+ (p'_{2})^{2})} \eea
is finite. Label the linear combinations of $p_i$'s and $p'_i$'s in
the denominator as $v_1, ... , v_5$. We see that \be
\lll{leg}(p_1,p'_1,p_2,p'_2) = (v_2 - v_3, v_3 - v_4, v_5 - v_4,v_5)
\ee Substituting these values into the integrand, we see that each
$v_j$ appears to a maximum power of $2/3$ in the numerator. This
implies that (\ref{wri}) is bounded by
\be \lll{fin} c \int_{\reals^{4}} (\prod_{j=1}^5
\frac{1}{1+|v_j|^{4/3}}) (\prod_{i=1}^2 dp_i dp'_i) \ee
We may
transform linearly from $(p_1,p'_1,p_2,p'_2)$ to
$(v_2,v_3,v_4,v_5)$, as is shown by (\ref{leg}).The resulting
integral is finite, as the power of each variable in the denominator
is greater than 1.

The general case may be handled in exactly the same way. That is,
given an ordering of $s_i$'s, $p_i$'s, and $p'_i$'s in (\ref{hoo}),
we may rewrite the expectation so that it factors as in (\ref{ywa}).
We substitute $(z_1, ... , z_{3n})$ for the differences of $s_i$'s,
$p_i$'s and $p'_i$'s, where $z_{3n}$ is defined to be $0$ to
simplify what follows. We use the bound (\ref{bund}), and arrive at
an expression of the form
\be \lll{wait} \int_{\reals^{2n}} \frac{(\prod_{i=1}^n |p_{i}|^{1/3}
|p'_{i}|^{1/3})(\prod_{i=1}^n dp_i dp'_i)}{\prod_{j=1}^{3n}
(1+|v_j|^2)} \ee
Each $p_i$ and $p'_i$ can be expressed as $v_j - v_{j+1}$ for some
$j \geq 1$. To see that this is true, note that when we rewrite the
expectation, as in step (\ref{ywa}), the only terms containing the
$u_i$ corresponding to a given $p_i$ will be

\be ... + v_j (B_a - B_{u_i}) + v_{j+1} (B_{u_i} - B_b) + ... \ee
where $a$ and $b$ denote the $s$, $u$, or $u'$ appearing immediately
before or after $u_i$ on the region to be integrated over. Comparing
this with the coefficient of $B_{u_i}$ in (\ref{wri}), we see that
$p_i = v_j - v_{j+1}$. Let us denote by $j(i)$ and $j'(i)$ as the
$j$ values for which $p_i = v_{j(i)} - v_{j(i)+1}$ and $p'_i =
v_{j'(i)} - v_{j'(i)+1}$. If we replace each $p_i$ and $p'_i$ in
(\ref{wait}) by the correct $v_j$, we see that (\ref{wait}) is
bounded by
\be \lll{wai} \int_{\reals^{2n}} \frac{(\prod_{i=1}^n
(|v_{j(i)}|+|v_{j(i)+1}|)^{1/3}(|v_{j'(i)}|+|v_{j'(i)+1}|)^{1/3})(\prod_{i=1}^n
dp_i dp'_i)}{\prod_{j=1}^{3n} (1+|v_j|^2)} \ee
Each $v_j$ appears at most twice in the numerator of (\ref{wai}), so
(\ref{wai}) is bounded by
\be \lll{wai2} \int_{\reals^{2n}}\frac{(\prod_{i=1}^n dp_i
dp'_i)}{\prod_{j=1}^{3n} (1+|v_j|^{4/3})} \ee This is finite, as the
set of $v_j$'s spans the set of $p_i$'s and $p'_i$'s.

This handles the variance in $x$ and $\eps$. We must still control
the variance in $t$. Assume $t'> t$. Then
\bea  && E(V(x,\eps,t)-V(x,\eps,t'))^{2n} \nn\\&& \leq
c\int_{\reals^{2n}}\int_{[t,t']^{n}} \int_{[o,s]^{2n}} E[exp(i\sum
_{i=1}^n[p_{i}(B_{s_{i}}-B_{u_{i}}) + p'_{i}(B_{s_{i}}-B_{u'_{i}})] \nn\\
\nn && (\prod_{i=1}^n du_{i} du'_{i} ds_{i} dp_{i} dp'_{i})\eea We
will follow the steps (\ref{wri}) through (\ref{fin}). Note however
that (\ref{bund}) may be combined with H\"{o}lder's inequality to
obtain
\be \lll{bund2} \int_{0}^{t'-t}e^{-rb^{2}}dr \leq
\frac{c|t-t'|^{1/p}}{(1+b^{2})^{1/q}} \ee
for any $p,q > 1$ such
that $1/p + 1/q = 1$. We will use (\ref{bund2}) in place of
(\ref{bund}), with $q = 10/9, p = 10$. Since all of the $s_i$'s must
be restricted to the interval $[t,t']$, we will have at least $n-1$
of the $z_k$'s restricted to $[0,t'-t]$(Recall that the $z_k$'s are
defined as in (\ref{jpb})). This shows that
\bea \lll{smp} \nn E(V(x,\eps,t)-V(x,\eps,t'))^{2n}  \\ \nn \leq
c|t-t'|^{(n-1)/10}\int_{\reals^{2n}} \frac{(\prod_{i=1}^n
|p_{i}|^{1/3}|p'_{i}|^{1/3})(\prod_{i=1}^n dp_i
dp'_i)}{\prod_{i=1}^{3n} (1+|v_j|^2)^{9/10}} \eea Following the
steps (\ref{wai}) and (\ref{wai2}), the integral in (\ref{smp}) is
bounded by

\be \lll{wai2} \int_{\reals^{2n}}\frac{(\prod_{i=1}^n dp_i
dp'_i)}{\prod_{j=1}^{3n} (1+|v_j|^{4/3-2/10})} \ee This integral is
finite, so (\ref{req2}) is proved. We have therefore proved

\be E[V(x,\eps,t)-V(x',\eps ',t')]^{2n} \leq C_M
|(x,\eps,t)-(x',\eps',t')|^{n/20}\ee By Kolmogorov's continuity
criterion, this implies that we may let $\eps \lar 0$ to obtain a
process which is defined on a set of full measure for all $x,t$.
That process has already been proved to be almost surely equal to
$\int_0^t L_s^{B_s - x}dB_s$ for each $x,t$. This completes the
proof.

\medskip

\noindent{\bf Acknowledgements}
\medskip

I am indebted to Jay Rosen for suggestion this problem to me, as
well as for all of his help and support throughout my graduate
career.


\begin{thebibliography}{1}

\bibitem{rosmar}
Marcus, M., Rosen, J.(2006) {\it Markov Processes, Gaussian
Processes, and Local Times}, preprint.

\bibitem{revyor}
Revuz, D., Yor, M. (1968) {\it Continuous Martingales and Brownian
Motion} Springer, Berlin.

\bibitem{rosen2}
Rosen, J. (2005). Derivatives of self-intersection local times, \,
{\it S\'{e}minaire de Probabilit\'{e}s,}\,\,XXXVIII,\,
Springer-Verlag, New York , LNM 1857, 171-184.

\end{thebibliography}

\def\noopsort#1{} \def\printfirst#1#2{#1} \def\singleletter#1{#1}
   \def\switchargs#1#2{#2#1} \def\bibsameauth{\leavevmode\vrule height .1ex
   depth 0pt width 2.3em\relax\,}
\makeatletter \renewcommand{\@biblabel}[1]{\hfill#1.}\makeatother

\it{Greg Markowsky}

\it{1 Edgewood Dr.}

\it{Orono, ME 04473 USA}

\it{greg@markowsky.com}

\end{document}